\documentclass[preprint,preprintnumbers,amsmath,amssymb]{revtex4}

\usepackage{graphicx}
\usepackage{dcolumn}
\usepackage{bm}

\newcommand{\Op}[2]{O_p(n^{- \frac{#1}{#2}} ) }
\newcommand{\On}[2]{O(n^{- \frac{#1}{#2} } ) }

\newcommand{\half}{\frac{1}{2} }
\newcommand{\real}{\mathbf{R}}

\newcommand{\Tr}{\mathrm{Tr}}

\newcommand{\rmd}{\mathrm{d}}

\newcommand{\Si}{\Sigma^{-1}}
\newcommand{\Sp}[1]{\frac{\partial \Sigma}{\partial \theta^{#1}}}
\newcommand{\Spp}[2]{\frac{\partial^2 \Sigma}{\partial \theta^{#1} \partial \theta^{#2}   }}
\newcommand{\Sppp}[3]{\frac{\partial^3 \Sigma}{\partial \theta^{#1}\partial \theta^{#2}\partial \theta^{#3} }}
\newcommand{\mGamma}{\stackrel{(m)}{\Gamma} }
\newcommand{\eGamma}{\stackrel{(e)}{\Gamma} } 

\newcommand{\mle}{\hat{\theta}}
\newcommand{\Smle}{\hat{S}}
\newcommand{\domega}{\frac{\mathrm{d} \omega}{4\pi} }

\newcommand{\Fvec}[1]{\left( \partial_{#1} \log \frac{f}{\pi_J} \right)}

\newcommand{\proof }{\textit{Proof.}\\}
\newcommand{\lemma}[1]{\textit{Lemma #1.}\\}
\newcommand{\proposition}[1]{\textit{Proposition #1.}\\}
\newcommand{\theorem}[1]{\textit{Theorem #1.}\\}
\newcommand{\qed}{\textit{Q.E.D.}\\}

\begin{document}

\preprint{APS/123-QED}

\title{Asymptotic Expansion of the Risk Difference \\ of the Bayesian Spectral Density \\ in the ARMA model}
     
\author{Fuyuhiko Tanaka}
\email{ftanaka@stat.t.u-tokyo.ac.jp}
\author{Fumiyasu Komaki}%
\affiliation{%
Department of Mathematical Informatics, University of Tokyo, 7-3-1 Hongo, Bunkyo-ku, Tokyo, 113-8656  Japan
}%

\date{\today}

\begin{abstract}
The autoregressive moving average (ARMA) model is one of the most important models in time series analysis.
We consider the Bayesian estimation of an unknown spectral density in the ARMA model.
In the i.i.d. cases, Komaki showed that Bayesian predictive densities based on a superharmonic prior
 asymptotically dominate those based on the Jeffreys prior~\cite{Komaki2002}.
It is shown by using the asymptotic expansion of the risk difference.
We obtain the corresponding result in the ARMA model.
\end{abstract}


\maketitle

\section{Introduction}

Let us consider a prediction problem in the Bayesian framework.
Suppose that a parametric model 
\[
 \mathcal{M}:= \{ p(x| \theta) : \theta \in \Theta \subseteq \real^k \}
\]
is given and our problem is to estimate $p(y| \theta)$ itself from an observation $x$.
If a proper prior density $\pi (\theta) $ is known, the best predictive density minimizing the average risk is obtained by~\cite{Ait1975}
\[
  p_{\pi }(y| x) := \int p(y|\theta) \pi(\theta | x ) \rmd \theta .
\]
If one has no knowledge on the unknown parameter $\theta $, he or she tends to adopt a noninformative prior.
It is often recommended to use the Jeffreys prior as a noninformative prior due to several reasons.
However  the Jeffreys prior is often improper, i.e., $\int \pi(\theta) \rmd \theta = \infty $.
In such a situation, the above result does not hold any more and other noninformative priors could be recommended.
One of those is a superharmonic prior.
Komaki showed that Bayesian predictive distributions based on a superharmonic prior asymptotically dominate those based on the Jeffreys prior~\cite{Komaki2002}.
He compared two prior distributions by using the asymptotic expansion of the risk difference.

In the present paper, we extend this result to the ARMA process.
We formulate the prediction problem of spectral densities in the ARMA model as described below and obtain the asymptotic expansion of the risk difference.
Our conclusion is the same as in the i.i.d. cases. 
Since we used the properties of the ARMA model only when evaluating the expectation of the log likelihood, 
it can be expected that almost all our arguments hold true in general stationary Gaussian processes. 

\subsection{General setting}

Let us consider a parametric model of stationary Gaussian process with mean zero.
It is known that a stationary Gaussian process corresponds to its spectral  density one-to-one (for proof, see, e.g., ~\cite{BD}).
Thus, we focus on the estimation of the true spectral density $S(\omega | \theta_0)$ in a parametric family of spectral densities
\[
 \mathcal{M}:= \{ S(\omega | \theta) : \theta \in \Theta \subseteq \real^k \}.
\] 
The performance of a spectral density estimator
 $\hat{S}(\omega) $ is evaluated by the Kullback-Leibler divergence. 
\[
 D( S(\omega | \theta_0) || \hat{S}(\omega) ) :=  \int _{-\pi}^{\pi} \domega \left\{   \frac{   S(\omega | \theta_0)       }{ \hat{S}(\omega)     }    -1
 - \log \left( \frac{S(\omega | \theta_0) }{ \hat{S}(\omega) } \right)  \right\}.
\]
The above setting is proposed by Komaki~\cite{Komaki1999}.

\subsection{Bayesian framework}

First, let us consider minimizing the average risk assuming that a proper prior density $\pi (\theta ) $ is known in advance.  
Aitchison's result~\cite{Ait1975} applies to this setting.
The spectral density estimator minimizing  the average risk,  
\begin{eqnarray*}
&&E^{\Theta} E^{X}[D(S(\omega | \theta ) || \hat{S}(\omega)  )  ]  \\
&& {} := \int \rmd \theta \pi (\theta)  \int \rmd x_1 \dots \rmd x_n  p_n(x_1 , \dots , x_n | \theta ) 
 D(S(\omega | \theta ) || \hat{S}(\omega ) ),
\end{eqnarray*}
 is given by the Bayesian spectral density (with respect to $\pi (\theta) $), which is defined by
\begin{equation}
 S_{\pi}(\omega) := \int S(\omega | \theta ) \pi(\theta | x )  \rmd \theta. \label{def:Bayesian}
\end{equation}
We call $S_{\pi}(\omega) $ in (\ref{def:Bayesian}) a \textit{Bayesian spectral density} even when an improper prior distribution is considered.

\subsection{Choice of a noninformative prior}

If  one has no information on the unknown parameter $\theta $, it is natural to adopt a noninformative prior in the Bayesian framework.
There is much room to argue the choice of a noninformative prior.
While the Jeffreys prior is a well-known candidate from several reasons, 
it can be expected that it is better to adopt a superharmonic prior in some cases.
The reason is that stationary Gaussian processes are getting  close to the  i.i.d. cases as the sample size becomes large
 and a superharmonic prior can be better than the Jeffreys prior in the i.i.d. cases.

\subsection{Construction}

In the following section, we briefly review basic results necessary to the asymptotic expansion.
The asymptotic expansion of the posterior distribution is presented.
In section 3, we obtain the asymptotic expansion of the Bayesian spectral density.
For the ARMA model, it can be also written in the differential-geometrical quantities as in the i.i.d. cases.
In section 4, we evaluate the expectation of the KL-divergence from the true spectral density $S(\omega | \theta_0)$ 
to the Bayesian spectral density $S_{f}$ up to the second order for an arbitrary prior (possibly improper) $f(\theta)$.
Finally, we obtain the principal term of the risk difference between $S_{f}$ and $S_{\pi_{J} }$, where $\pi_{J}$ denotes the Jeffreys prior.
As a direct consequence of this result, a superharmonic prior is recommended as a noninformative one if there exists a positive superharmonic function 
on the corresponding model manifold.

\section{Preliminary}

\subsection{Notation and assumption}

In the present section, we consider general stationary Gaussian processes with mean zero.
We recall that the likelihood function is given by
\begin{equation}
 p_n(x_1, \cdots ,  x_n| \theta ) = \frac{1}{ \sqrt{\det\Sigma_n(\theta) (2\pi)^n  }  } \exp \left( - \frac{1}{2} X'_n \Sigma_n(\theta)^{-1} X_n  \right),  \label{eq:dist}
\end{equation}
where $X_n = (x_1, \cdots, x_n)$ and $\Sigma_n$ denotes a covariance matrix and
$\theta \in \Theta \subset \real^k$ denotes an unknown parameter.
We often use the log likelihood of the form omitting the constant term
\[
 l_n (\theta ) = - \frac{1}{2} X'_n \Sigma_n^{-1}(\theta ) X_n - \frac{1}{2} \log \det \Sigma_n(\theta). 
\]
We also assume that an arbitrary prior density $\pi(\theta)$ on $\Theta$ is given.\\
\\
\\
In the present paper, we assume several regularity conditions (See, e.g., Taniguchi and Kakizawa~\cite{Taniguchi}).
\\
\\
Differential operators are denoted as $\partial_j := \frac{\partial }{\partial \theta^j}$ as often seen in the differential geometry
 (for other basic notation,  see, e.g., Kobayashi and Nomizu~\cite{KN}).
We also use Einstein's summation convention: if an index occurs twice in any one term, once as an upper and once as a lower index, 
summation over that index is implied.

\subsection{Asymptotic expansion of the posterior density}

For stationary Gaussian time series models, we have the asymptotic expansion of the posterior density
\[
 \pi(\theta | X_n) = \frac{ p_n(X_n | \theta) \pi(\theta) }{ \int \rmd \theta p_n(X_n | \theta) \pi(\theta)    }.
\]
Here, we give the moment form of the expansion.\\
\\
\lemma{1}
The asymptotic expansion of the posterior density in the moment form is given by
\begin{eqnarray}
&& E^{\pi}[ \Delta \theta ^{i_1} \cdots \Delta \theta^{i_p}  ]   \nonumber    \\
&&:= \int  \Delta \theta ^{i_1} \cdots \Delta \theta^{i_p} 
 \pi(\theta | x) \rmd \theta \nonumber     \\
			&&= \int   \Delta \theta^{i_1} \cdots \Delta \theta^{i_p} 
\frac{1}{\sqrt{\det \left( n^{-1} J_{n}(\mle)^{-1} \right)   } } 
			\frac{1}{(2\pi)^{k/2} } 
			\exp \left( - \frac{n}{2}J_n(\mle)(\theta - \mle)^2  \right) \nonumber \\			
			&& {} \times \left\{ 1 + \frac{A_n(x)}{\sqrt{n}} + O_p(n^{-1})  \right\}  \rmd \theta,   \nonumber \\ 
 && \label{eq:lemma1}
\end{eqnarray}
where $\mle $ is the maximum likelihood estimate and
\[
\Delta \theta^i := \theta^i - \mle^i (x), \  
 [J_{n}(\mle)]_{ij} := - \frac{1}{n} \frac{ \partial^2 l_n (\mle ) }{ \partial \theta^i \partial \theta^j}
\] 
\[
 A_n(x) := \frac{1}{3!} \frac{\partial_{i_1} \partial_{i_2} \partial_{i_3} l_n(\mle )}{n} ( \sqrt{n} \Delta \theta^{i_1} )( \sqrt{n} \Delta \theta^{i_2} )
( \sqrt{n} \Delta \theta^{i_3} )
	+ \partial_j \log \pi (\mle) (\sqrt{n} \Delta \theta^{j}    ) 
\]
For derivation, see, for example, Philippe and Rousseau~\cite{PandR}.
From this formula, we can calculate $E^{\pi}[ \Delta \theta^{i_1} \cdots  \Delta \theta^{i_p}  ]$ in an arbitrary order.
The $p$-th moment is defined by 
\begin{eqnarray*}
 I^{i_1 i_2 \dots i_p}(\mle ) &:=& \int \cdots \int y^{i_1} y^{i_2} \cdots y^{i_p} 
  \frac{1}{\sqrt{\det \left( n^{-1} J_{n}(\mle)^{-1} \right)   } } 
			\frac{1}{(2\pi)^{k/2} } \\
			& & \times \exp \left( - \frac{n}{2}J_n(\mle)_{lm} y^l y^m  \right) \rmd y^{1}\dots \rmd y^{k}
\end{eqnarray*}
We need only the second moment $I^{ij}(\mle )$ and the fourth moment $I^{ijkl}(\mle) $ in the present paper.
\[
 I^{i_1 i_2} (\mle ) = \left( \frac{J_n^{-1}( \mle) }{n} \right)_{i_1 i_2}
\]
\[
 I^{i_1 i_2 i_3 i_4}(\mle) =I^{i_1 i_2} (\mle ) I^{i_3 i_4} (\mle )
				+ I^{i_1 i_3} (\mle ) I^{i_2 i_4} (\mle )
				 + I^{i_1 i_4} (\mle ) I^{i_2 i_3} (\mle )
\]
Stochastic order of them is evaluated as (for even $p$)
\[
 I^{i_1 i_2} = O_p(n^{-1}),   \quad I^{i_1 i_2 i_3 i_4 } = O_p(n^{-2}) , \cdots , I^{i_1 \cdots i_p} =  O_p(n^{- \frac{p}{2}}) .
\]
Using these moment formula, we obtain
\[
\begin{array}{lcl}
 E^{\pi}[\Delta \theta^i] \equiv B^i _{\pi}(\mle) &=& \frac{1}{ \sqrt{n}} \frac{1}{3!} \left( \frac{  \partial_{i_1} \partial_{i_2}\partial_{i_3}   l_n(\mle ) }{n} \right)
    (\sqrt{n})^3 I^{i_1 i_2 i_3 i}(\mle)
 \\
&& + \frac{1}{\sqrt{n}}  \partial_{i_1} \log \pi  (\mle) \sqrt{n} I^{i_1 i}(\mle) + O_p(n^{-2})  (= O_{p}(n^{-1})  ) \\
E^{\pi}[\Delta \theta^i \Delta \theta^j ] &=& I^{ij}(\mle) + O_p(n^{-2}) \\
E^{\pi}[\Delta \theta^i \Delta \theta^j \Delta \theta^k ] &=& O_p (n^{-2}) \\
\end{array}
\]

\section{Asymptotic Expansion of the Spectral Density}

We consider a parametric family of spectral density 
\[
\mathcal{M} := \{ S(\omega | \theta) : \theta \in \Theta \subseteq \real^k, S \mbox{ corresponding to an ARMA process } \}.
\]
Let  $f (\theta) $ be an arbitrary prior distribution on $\Theta $ and $\pi_{J}(\theta ) $ be the Jeffreys prior.
From now on, $\theta_0 $ denotes the true parameter. 
We consider estimating the true spectral density $S(\omega | \theta_0)$ itself instead of $\theta_0 $.
From the $n$ data, $X_n := (x_1, \dots, x_n )$ subject to an ARMA process, we construct the posterior distribution $f(\theta |X_n)$ and
the Bayesian spectral density $\hat{S}_{f}(\omega) $ with respect to $f(\theta)$ is given by 
\[
 \hat{S}_{f}(\omega) := \int  S(\omega | \theta ) f(\theta | X_n)  \rmd \theta .
\]
From the result of the previous section, we obtain\\
\\
\lemma{2}
Let the maximum likelihood estimator $\mle = \theta_0 + \Op{1}{2}$ given, then the Bayesian spectral density is evaluated as
\begin{equation}
 \hat{S}_{f}( \omega) = S(\omega | \mle ) + \partial_i S(\omega | \mle) B^{i}_{f}(\mle)
	+ \half \partial_i \partial_j S(\omega | \mle ) I^{ij}(\mle) + \Op{3}{2}.  \label{eq:lemma2}
\end{equation}
\\
\proof 
Using the Taylor expansion of $S(\omega | \theta )$ around $\mle $, Eq.(\ref{eq:lemma2}) immeadiately follows from Lemma1.\\
\qed \\

\subsection{Geometrical expansion of the derivatives of the log likelihood}

We rewrite the Eq.(\ref{eq:lemma2}) using geometrical quantities, which are defined by 
\[
\begin{array}{lcl}
 g_{ij} &=& \int \domega  \partial_i \log S(\omega |\theta_0)  \partial_j \log S(\omega |\theta_0) \\
 \mGamma_{i,jk} &=&  \int \domega  \partial_i \log S(\omega |\theta_0) 
   \frac{ \partial_j \partial_k S(\omega |\theta_0)  }{S(\omega |\theta_0) } \\
 M_{i,jkl} &=&  \int \domega  \partial_i \log S(\omega |\theta_0)  \frac{ \partial_j \partial_k \partial_l S(\omega |\theta_0) }{S(\omega |\theta_0) } \\
 N_{ij,kl} &=&  \int \domega \frac{ \partial_i \partial_j S(\omega |\theta_0)  }{S(\omega |\theta_0) }
\frac{ \partial_k \partial_l  S(\omega |\theta_0) }{S(\omega |\theta_0) }\\
    T_{ijk} &=& \int \frac{ \rmd \omega}{2\pi }  \partial_i \log S(\omega |\theta_0) \partial_j \log S(\omega |\theta_0) \partial_k \log S(\omega |\theta_0) \\
   L_{ij,kl} &=& \int \domega \partial_i \log S(\omega |\theta_0) \partial_j \log S(\omega |\theta_0)  \frac{ \partial_k \partial_l S(\omega |\theta_0) }{S(\omega |\theta_0) } \\
\end{array}
\]
For these geometrical notations, e.g., see Amari~\cite{Amari1987, Amari2000}.
It is convenient to introduce some notation for the log likelihood $l_n (\theta_0) $ and its derivatives.
\[
 L_{i_1 \cdots i_p}(\theta_0) := \frac{1}{n} \frac{\partial^p l_n(\theta_0) }{\partial \theta^{i_1} \cdots \partial \theta^{i_p} } = O_p(1)
\]
and
\[
 m_{i_1 \cdots i_p}(\theta ) := E_{\theta_0} [L_{i_1\cdots i_p}(\theta )]
\]
Note that $\theta \neq \theta_0 $, for example, $m_i (\theta )= E_{\theta_0}[L_i (\theta)] \neq 0$ but $m_i (\theta_0) = 0$.
Likewise, the expectation of the product of the log derivatives are defined by
\[
 m_{i_1 \cdots i_p, j_1 \cdots j_q}(\theta) := E_{\theta_0}[ L_{i_1\cdots i_p}(\theta ) L_{j_1\cdots j_q}(\theta )].
\]
We omit the argument $\theta_0$ if otherwise necessary.
Other important notations are $L^{ij}(\theta )$ and $m^{ij}(\theta)$.
Each of them denotes the inverse matrix of $L_{ij}(\theta)$  and that of $m_{ij}(\theta) $.
%
%
%
%
%
%
%
%
Note that $
L^{ij} := (L^{-1})_{ij} = m^{ij} - m^{il} (\delta L)_{lk} m^{kj} + \cdots,  
$
where $(\delta L)_{lk} = L_{lk} - m_{lk}(= \Op{1}{2})$.\\
\\
\lemma{3}
For the ARMA model, we obtain the explicit forms of $m_{ij},m^{ij}, m_{ijk}$ and $m_{ij,k}$. They are represented by geometrical quantities. 
\[
 \begin{array}{lclcl}
 m_{ij} &=&  -g_{ij} +O(n^{-1})  &=& O(1)\\
 m^{ij} &=&  -g^{ij} +O(n^{-1}) &=& O(1)\\
 m_{ijk} &=&   2T_{ijk} - ( \mGamma_{i,jk} + \mGamma_{j,ik} + \mGamma_{k,ij} ) + O(n^{-1}) & &  \\
         &=&  - ( \eGamma_{i,jk} + \eGamma_{j,ik} + \eGamma_{k,ij} + T_{ijk}) + O(n^{-1}) &=&  O(1) \\
 n m_{ij,k} &=& \mGamma_{k,ij} - T_{ijk} +O(n^{-1}) = \eGamma_{k,ij}+O(n^{-1}) &=& O(1) \\
\end{array}
\]
where $\eGamma_{i,jk} := \mGamma_{i,jk} - T_{ijk}$.\\
\\
\proof
First, we show $m_{ij} =- g_{ij} + O(n^{-1})$. 
From straightforward calculation, we obtain
\begin{equation}
m_{ij} = - \frac{1}{2n}  \Tr \left(\Si \Sp{i} \Si \Sp{j}\right) \label{eq:mij}
\end{equation}
Here, the following fact holds for a parametric family of the spectral density of the ARMA model (See, Lemma 4.1.2~\cite{Taniguchi}), \\
\\
\textit{Fact.}~\cite{Taniguchi}\\
Let $p_n(x_1, \dots, x_n| \theta) $ is given by Eq.(\ref{eq:dist}) and the corresponding spectral density is $S(\omega | \theta)$.
Then,  for arbitrary $k \geq 1$ and $p_1, \dots , p_k$, 
\begin{eqnarray*}
 &&\frac{1}{n}\Tr  \left \{ \! \Sigma^{-1} (D_{p_1} \Sigma ) \Sigma^{-1} \cdots (D_{p_k} \Sigma) \!  \right\} 
  \\
 &=& \frac{1}{2\pi}\!\! \int_{-\pi}^{\pi} \! \frac{ D_{p_1}S(\omega | \theta) }{S(\omega | \theta)} \cdots  
\frac{  D_{p_k}S(\omega | \theta) }{ S(\omega | \theta) }  \rmd \omega + O(n^{-1}),
\end{eqnarray*}
where $D_{p} $ denotes an arbitrary $p$-th order differential operator $\frac{\partial}{\partial \theta^{l_1} }\cdots \frac{\partial}{\partial \theta^{l_p}}$.\\
\\
Using the fact, the trace in the r.h.s. of (\ref{eq:mij}) is rewritten in the form of the integral  
\[
 \frac{1}{n} \Tr \left( \Si \Sp{i} \Si  \Sp{j} \right) 
   = \frac{1}{2\pi} \int_{-\pi}^{\pi} \frac{ \partial_i S(\omega | \theta) }{S(\omega | \theta) }\frac{ \partial_j S(\omega | \theta) }{S(\omega | \theta) }
 \mathrm{d} \omega + O(n^{-1}).
\]
Thus, we obtain
\[
 m_{ij} = - g_{ij} + O(n^{-1}).
\]
Since $m^{ij}$ is the inverse of  $m_{ij}$, the second equation clearly holds.
The other equations are shown in the same way.\\
\qed

\subsection{Geometrical expansion of  the Bayesian spectral density}

Now we rewrite the asymptotic expansion (\ref{eq:lemma2}) in the geometrical quantities.\\

\noindent\lemma{4}
Let $\mle $ be the maximum likelihood estimate, then the following expansion holds.
\begin{eqnarray}
 \hat{S}_f (\omega) &=& S(\omega | \mle) + \frac{1}{2n} g^{ij}(\mle) \left( \partial_i \partial_j S(\omega | \mle) - \mGamma_{ij}^{k}(\mle) \partial_k S( \omega | \mle)  \right) \nonumber  \\
		&& + \frac{1}{n} g^{ij}(\mle) \left\{  \partial_i \log \frac{f}{\pi_J}(\mle) + \half T_i (\mle)   \right\} \partial_j S(\omega | \mle) + \Op{3}{2},\label{eq:lemma4}
\end{eqnarray}
where $T_i := T_{ijk} g^{jk}$.
Note that Eq.(\ref{eq:lemma4}) is formally in the same form as those in the i.i.d. cases  if one reads $p(y |\theta)$ as $S(\omega | \theta)$.(See, Komaki\cite{Komaki2002}).  \\
\\
\proof
In order to prove Eq.(\ref{eq:lemma4}), one can neglect $\Op{3}{2}$ terms.
For example, up to this order, the following identity holds
\begin{eqnarray*}
 I^{ij}(\mle) = -\frac{1}{n} L^{ij} (\mle)& =& -\frac{1}{n} m^{ij}(\theta_0 ) + \Op{3}{2} \\
&= &\frac{1}{n} g^{ij}(\theta_0) + \Op{3}{2} \\
& =& \frac{1}{n} g^{ij} (\mle) + \Op{3}{2}.
\end{eqnarray*}
Now let us rewrite the principal term of  $\partial_{i} S(\omega | \mle) B^{i}_{f} (\mle)$ in Eq.(\ref{eq:lemma2}).
From Lemma 3,
\begin{eqnarray*}
L_{jkl}(\mle) &=& m_{jkl}(\theta_0) + \Op{1}{2} 
\\
		&=& \left\{  2T_{jkl}(\theta_0 ) - (\mGamma_{j,kl}(\theta_0 ) + \mGamma_{k,jl} (\theta_0 )  + \mGamma_{l,jk} (\theta_0 )  )  \right\} + \Op{1}{2} \\
		&=& \left\{  2T_{jkl}(\mle) - (\mGamma_{j,kl}(\mle ) + \mGamma_{k,jl} (\mle )  + \mGamma_{l,jk} (\mle )  )  \right\} + \Op{1}{2} 
\end{eqnarray*}
and
\begin{eqnarray*}
&&L_{jkl}(\mle) L^{ij}(\mle) L^{kl}(\mle) \\
&=& \left\{  2T_{jkl}(\mle) - (\mGamma_{j,kl}(\mle ) + \mGamma_{k,jl} (\mle )  + \mGamma_{l,jk} (\mle)  )  \right\} 
    g^{ij} (\mle) g^{kl}(\mle) + \Op{1}{2} \\
		&=& - \mGamma_{kl}^{i}(\mle) g^{kl}(\mle) - ( \eGamma_{k,jl}(\mle) +\eGamma_{l,jk}(\mle)) g^{kl}(\mle) g^{ij}(\mle)  + \Op{1}{2} \\
		&=& - \mGamma_{kl}^{i}(\mle) g^{kl}(\mle)  - 2 \eGamma_{jl}^{l} (\mle) g^{ij} (\mle) +  \Op{1}{2} .
\end{eqnarray*}
Thus, 
\begin{eqnarray}
nB_{f}^{i}(\mle) \partial_i S( \omega |  \mle) &=& \half \partial_i S(\omega |  \mle) \left\{  - \mGamma_{jk}^{i}(\mle) g^{jk}(\mle) - 2 \eGamma_{jk}^{k} (\mle) g^{ij} (\mle) +  \Op{1}{2} 
  \right\} \nonumber \\
		&&  + g^{ij}(\mle)  \partial_j \log f (\mle) \partial_i S(\omega |  \mle)  + \Op{1}{2} \nonumber \\    
		&=& -\half g^{jk} (\mle) \mGamma_{jk}^{i} (\mle) \partial_i S(\omega |  \mle)  \nonumber \\
		&  &  + g^{ij}(\mle)  \partial_i S(\omega |  \mle )  \left\{  \partial_j \log f (\mle) - \eGamma_{jk}^{k} (\mle)  \right\}   + \Op{1}{2} \nonumber \\    
		&=& -\half g^{jk} (\mle) \mGamma_{jk}^{i} (\mle) \partial_i S(\omega |  \mle)  \nonumber \\
		&  &  + g^{ij}(\mle)  \partial_i S(\omega |  \mle)  \left\{  \partial_j \log \frac{f}{\pi_{J}}(\mle) + \half T_j (\mle)  \right\}   + \Op{1}{2}. \label{eq:proof3} 
\end{eqnarray}
In the last equality, we used the relation
\[
\partial_i \log \pi_J = \partial_i \log \sqrt{g} = \Gamma_{ij}^{j} = \eGamma_{ij}^{j}+ \half T_{i}.
\]
Substituting Eq.(\ref{eq:proof3}) into Eq.(\ref{eq:lemma2}), we obtain
\begin{eqnarray*}
 \hat{S}_f (\omega) &=& S(\omega | \mle) + \frac{1}{2n} g^{ij}(\mle) \left( \partial_i \partial_j S(\omega | \mle) - \mGamma_{ij}^{k}(\mle) \partial_k S( \omega | \mle)  \right) \\
		&& + \frac{1}{n} g^{ij}(\mle) \left\{  \partial_i \log \frac{f}{\pi_J}(\mle) + \half T_i (\mle)   \right\} \partial_j S(\omega | \mle) + \Op{3}{2}.
\end{eqnarray*}
\qed

\section{Asymptotic Expansion of the Expectation of the KL-divergence}

In this section, we evaluate the expectation of the KL-divergence up to the second order, focusing on the terms including the prior distribution $f(\theta)$.
For simplicity, we introduce the following notation.
\[
S_0 := S(\omega | \theta_0), \ \Smle := S(\omega | \mle), \ S_f := \hat{S}_f (\omega )  
\]
and 
\[
 \frac{S_f - S_0}{ S_0} = \left( \frac{S_f  - \Smle }{S_0} \right) +  \left( \frac{\Smle  - S_0}{S_0} \right) =: \Delta S_f + \Delta S_m.
\]
While the first term depends on the prior $f$ the second one is independent of  $f$.
Note that $\Delta S_m := \frac{\Smle -  S_0 }{S_0}=  \Op{1}{2}$ and $\Delta S_f := \frac{S_f - \Smle  }{S_0} = O_p(n^{-1}) $.
The KL-divergence from the true spectral density $S(\omega | \theta_0) $ to a Bayesian spectral density $\hat{S}_{f}(\omega)   $ is given by
\begin{eqnarray*}
 D( S_0 ||S_{f} ) &=&  \int _{-\pi}^{\pi} \domega \left\{   \frac{   S_0       }{ S_{f}     }    -1
 - \log \left( \frac{S_0 }{ S_{f}} \right)  \right\}  \\
&=& \int _{-\pi}^{\pi} \domega \left\{  \frac{1}{1+y} - 1- \log \left( \frac{1}{  1+y } \right)  \right\},
\end{eqnarray*}
where $y =  \Delta S_m + \Delta S_f $.
Due to the Taylor expansion,
\[
 \frac{1}{1+y} - 1- \log \left( \frac{1}{  1+y } \right) =  \frac{1}{2} y^2 - \frac{2}{3} y^3 + \frac{3}{4} y^4 + \cdots
= \sum_{n=2}^{\infty} (-1)^{n} \frac{n-1}{n} y^n ,
\]
we obtain
\begin{eqnarray*}
D(S_0 || S_f) &=& \half \left\{ 
		\int (\Delta S_m )^2 \domega + 2 \int \Delta S_m \Delta S_f \domega + \int (\Delta S_f )^2 \domega 
\right\} \\
	& & -\frac{2}{3} \left\{ \int (\Delta S_m )^3  \domega + 3 \int (\Delta S_m )^2 (\Delta S_f)  \right\} \\
	& & +\frac{3}{4} \left\{ \int (\Delta S_m)^4 \domega \right\} + \Op{5}{2} \\
	&=& U + \half V - 2 W \\
	&&    {} + \mbox{the terms independent of $f$} + \Op{5}{2} .
\end{eqnarray*}
We consider the expectation of  the following three terms including $f$:
\[
\begin{array}{lll}
 U &=& \int \Delta S_m \Delta S_f \domega ( = \Op{3}{2}) , \\
 V &=& \int (\Delta S_f )^2 \domega ( = O_p (n^{-2}) ) ,\\
 W &=& \int (\Delta S_m )^2 (\Delta S_f) \domega  ( = O_p(n^{-2}) ). \\
\end{array}
\]
Both $\Delta S_m$ and $\Delta S_f$ are given by 
\[
\Delta S_m = \frac{\partial_i S_0}{ S_0} \delta^i + \frac{1}{2} \left( \frac{\partial_i \partial_j S_0 }{S_0 } \delta^i \delta^j  \right)
 + \Op{3}{2}  
\]
and 
\[
\Delta S_f = \frac{1}{2n} \frac{1}{S_0} g^{ij}(\mle) ( \partial_i \partial_j \hat{S} - \mGamma_{ij}^{k} (\mle) \partial_k \hat{S} ) 
 + \frac{1}{n} g^{ij}(\mle) \frac{\partial_i \hat{S} }{S_0} F_j(\mle)  + \Op{3}{2},
\]
where
\[
\delta^i := \mle^i - \theta_0^i \ \mbox{and} \ F_j (\theta ) := \partial_j \log  \frac{f}{ \pi_J}  (\theta) + \half T_j(\theta ).
\]
It is convenient to use some formula for the maximum likelihood estimate $\mle$.
\[
\begin{array}{lcl}
\delta^i \delta^j &=& E_{\theta_0}[\delta^i \delta^j ] + \Op{3}{2} = \frac{1}{n} g^{ij}(\theta_0) + \Op{3}{2}, \\
E_{\theta_0}[\delta^i] &=& -\frac{1}{2n} \mGamma_{jk}^{i}(\theta_0) g^{jk}(\theta_0) + \On{3}{2}. \\ 
\end{array}
\]

\subsection{Evaluation of $V$ and $W$}

First of all, we calculate $V$ and $W$. We need the principal terms of each.
In the principal order, one can replace $\mle $ with $\theta_0$.
Each quantity is evaluated at the point $\theta_0$, but for simplicity, we omit $\theta_0$. 
\begin{eqnarray*}
n^2 V &=& \int (n\Delta S_f )^2 \domega \\
	&=& \int \domega  \left\{ \half g^{ij} \left( \frac{\partial_i \partial_j S}{S} - \mGamma_{ij}^{k} \frac{\partial_k S}{S} \right) + g^{ij}\frac{\partial_j S}{S} F_i \right\} \\
	&& \quad \quad \times \left\{ \half g^{ij} \left( \frac{\partial_i \partial_j S}{S} - \mGamma_{ij}^{k} \frac{\partial_k S}{S} \right) + g^{ij}\frac{\partial_j S}{S} F_i 
    \right\}  \\
	&=& g^{ij} F_i F_j + \frac{1}{4}g^{ij} g^{kl} N_{ij,kl} - \frac{1}{4} (\mGamma_{ij}^{l}g^{ij }) (\mGamma_{i'j'}^{l'}g^{i'j' })g_{ll'} + \Op{1}{2}.
\end{eqnarray*}
In the same manner, 
\begin{eqnarray*}
n^2 W &=& \int n (\Delta S_m)^2 (n\Delta S_f ) \domega \\
	&=& \int \domega \left\{ 
		n \delta^i \delta^j  \left( \frac{\partial_i S}{S}\right)\left( \frac{\partial_j S}{S}\right) + \Op{1}{2} \right\} \\
	& & \quad \quad  \times \left\{ 
		\half g^{kl} \left( \frac{\partial_k \partial_l S}{S} - \mGamma_{kl}^{m} \frac{\partial_m S }{S} + g^{kl} \frac{\partial_kS}{S} F_l  + \Op{1}{2} \right)
\right\}	\\
	&=& (n \delta^i \delta^j ) \left\{  
		\half g^{kl} (L_{ij,kl} - \mGamma_{kl}^{m}\half T_{ijm}  ) + g^{kl} F_l \half T_{ijk} 
\right\} + \Op{1}{2} \\
	&=& \half g^{ij} g^{kl} (L_{ij,kl} - \mGamma_{kl}^{m}\half T_{ijm}  )  + \half g^{ij} g^{kl} F_l T_{ijk} +\Op{1}{2}.
\end{eqnarray*}

\subsection{Evaluation of $U$}

For $U$, we need to evaluate terms up to the second principal order.
However, when it comes to the expectation, it is not so difficult.
For simplicity, we set
\[
A :=  \half g^{kl} \left( \partial_k \partial_l S - \mGamma_{kl}^{m} \partial_m S \right)   + g^{kl} \partial_k S \ F_l.
\]
We decompose $U$ into the following three terms:
\begin{eqnarray*}
U &=& \int \domega \frac{\partial_i S}{S} \frac{A}{S} \delta^i  \ (=: U_1) \\
	& & \quad +\half \delta^i \delta^j \int \domega \left( \frac{\partial_i \partial_j S}{S} \frac{A}{S}  \right) \ (=: U_2 )\\ 
	& & \quad + \delta^i \delta^j \int \domega \left( \frac{\partial_i S}{S} \frac{\partial_j A}{S}  \right) \  (=: U_3)
\end{eqnarray*}

\subsubsection{Evaluation of  $U_2$ and $U_3$}

Straightforward calculation yields 
\begin{eqnarray*}
U_2 &=&  \half \delta^i \delta^j \int \domega \left( \frac{\partial_i \partial_j S}{S} \frac{A}{S}  \right) \\ 
	&=& \frac{1}{n^2} \left\{  \frac{1}{4}g^{ij}g^{kl} N_{ij, kl} - \frac{1}{4}g^{ij} \mGamma_{k,ij} \mGamma_{pq}^{k} g^{pq} \right. \\
&& \quad \quad + \left. \half g^{ij} g^{kl}\mGamma_{l,ij} F_k \right\} + \Op{5}{2}
\end{eqnarray*}
and
\begin{eqnarray*}
U_3 &=& \delta^i \delta^j \int \domega \left( \frac{\partial_i S_0}{S_0} \frac{\partial_jA_0}{S_0}  \right) \\
	&=& \left( \frac{1}{n}g^{ij} + \Op{3}{2} \right) \left\{ 
		\frac{1}{2n} \left( g^{lk} M_{i, ljk} + \partial_j g^{lk} \mGamma_{i,lk} \right) 
 \right. \\
	& & - \frac{1}{2n} \left( 
	\mGamma_{i,jk} \mGamma_{pq}^{k} g^{pq} + g_{ik} \partial_{j} (g^{pq} \mGamma_{pq}^{k} ) 
\right) \\
	&& +\left. \frac{1}{n} \left( \mGamma_{i,jk} g^{kl}F_l + g_{ik} \partial_j (g^{kl} F_l)  \right)  + \Op{3}{2}\right\} \\
	&=& \frac{1}{2n^2} g^{ij}\left( g^{lk} M_{i, ljk} + \partial_j g^{lk} \mGamma_{i,lk}  - \mGamma_{i,jk} \mGamma_{pq}^{k} g^{pq} + g_{ik} \partial_{j} (g^{pq} \mGamma_{pq}^{k} ) 
 \right)  \\
	&& + \frac{1}{n^2} \left( g^{ij} g^{kl} \mGamma_{i,jk} F_l + \partial_k (g^{kl} F_l)  \right)  + \Op{5}{2}.
\end{eqnarray*}

\subsubsection{Evaluation of $U_1$}

Finally, we deal with $U_1$.

\begin{eqnarray*}
U_1 &=&\left( \int \domega \frac{\partial_i S}{S} \frac{A}{S} \right) \delta^i \\
 &=& \left( \frac{1}{n} F_i + T_i + O_p(n^{-2}) \right) \delta^i,
\end{eqnarray*}
where $T_i$ denotes $\Op{3}{2}$ terms.
Thus, the stochastic expansion of the $U_1$ term requires the second order of the asymptotic expansion of the posterior density.
However the evaluation of the expectation requires no such higher order terms because $E[T_i \delta^i] = \On{5}{2}$ in spite of  $T_i \delta^i = O_p(n^{-2})$.
Indeed, 
\begin{eqnarray*}
E[\delta^i T_i] &=& E[(\delta^i - E[\delta^i] ) \times T_i ] + E[\delta^i ] \times E[T_i] \\
		&=&  O (n^{-1}) \cdot \On{3}{2}  + O(n^{-1}) \cdot \On{3}{2} \\
		 &=& \On{5}{2}. 
\end{eqnarray*}
Thus, we obtain
\begin{eqnarray*}
E[U_1] &=& \frac{1}{n} F_i E[\delta^i] + E[T_i \delta^i ]\On{5}{2} \\
	&=& \frac{1}{n} F_i \left( - \frac{1}{2n} g^{im} g^{kj} \mGamma_{m,kj}  \right) + \On{5}{2} \\
	&=& - \frac{1}{2n^2} F_i \mGamma_{kj}^{i} g^{kj} + \On{5}{2}. 
\end{eqnarray*}

\subsection{Asymptotic expansion of the expectation of the KL-divergence}

Collecting the whole terms $U$,$V$ and $W$, we obtain
\begin{eqnarray*}
&& E_{\theta_0}[D(S_0 || S_f) ] \\
&=& - \frac{1}{2n^2} F_i \mGamma_{kj}^{i}g^{kj} + \frac{1}{2n^2} g^{ij} \mGamma^{k}_{ij} F_k+ \frac{1}{n^2} \left\{  g^{kl} \mGamma_{jk}^{j}F_l + \partial_k (g^{kl} F_l)  \right\} \\
&& + \frac{1}{2n^2} g^{ij} F_i F_j - \frac{1}{n^2} g^{kl} F_l T_k +\mbox{the terms independent of $f$} + \On{5}{2} \\
&=& \frac{1}{2n^2} g^{ij} F_i F_j  + \frac{1}{n^2} \eGamma_{lk}^{l} g^{kj}F_j + \frac{1}{n^2} \partial_k (g^{kj} F_j) \\
&& \quad \quad +\mbox{the terms independent of $f$} + \On{5}{2}\\
&=& \frac{1}{2n^2} g^{ij} F_i F_j + \frac{1}{n^2} \stackrel{e}{\nabla}_k(g^{kj} F_j) +\mbox{the terms independent of $f$} + \On{5}{2}.
\end{eqnarray*}
Summarizing this, we obtain the following proposition.\\
\\
\proposition{1}
Let $S_0 $ a true spectral density and $S_f$ the Bayesian spectral density with respect to $f(\theta)$. 
Then, the asymptotic expansion of the expectation of the KL-divergence from $S_0$ to $S_f$ is given by
\begin{eqnarray*}
 E_{\theta_0}[D(S_0 || S_f) ] 
 & =& \frac{1}{2n^2} g^{ij} 
 \left( \partial_i \log  \frac{f}{ \pi_J} + \half T_i \right) 
\left(  \partial_j \log  \frac{f}{ \pi_J} + \half T_j \right) \\
 && 	{}+ \frac{1}{n^2} \stackrel{e}{\nabla}_k  \left( 
		g^{kj}  \partial_j \log  \frac{f}{ \pi_J}  + \half T_j  \right) \\
 &&  \quad \quad +\mbox{the terms independent of $f$} + \On{5}{2}.
\end{eqnarray*}

\section{Comparison between $S_{\pi_{J}} $ and  $S_{f}$}

From the result in the previous section, we obtain the same result as that in the i.i.d. cases, which was shown by Komaki\cite{Komaki2002}.
Let us calculate the risk difference between two Bayesian spectral densities, one of which is based on the Jeffreys prior $\pi_{J}(\theta)$ and
 the other is based on an arbitrary prior $f$, 
\begin{eqnarray*}
&&n^2 E_{\theta_0} [ D(S_0 || \hat{S}_{\pi_{J}} )  ] - n^2 E_{\theta_0} [ D(S_0 || \hat{S}_{f} )  ]  \\
&=&  \frac{1}{8} g^{ij} T_i T_j + \stackrel{e}{\nabla}_k(g^{kj} \half T_j) - \frac{1}{2} g^{ij} \left( \partial_i \log \frac{f}{\pi_J} + \half T_i \right) \left( \partial_j \log \frac{f}{\pi_J} +\half T_j \right)  \\
&&    -  \stackrel{e}{\nabla}_k \left\{ g^{kj}  \left( \partial_j \log \frac{f}{\pi_J} + \half T_j \right)   \right\} + \On{1}{2}  \\
&=& - \frac{1}{2} g^{ij} \left( \partial_i \log \frac{f}{\pi_J} \right) \left( \partial_j \log \frac{f}{\pi_J} \right) - \frac{1}{2} g^{ij} T_i \partial_j \log \frac{f}{\pi_J} \\
&&	-  \stackrel{e}{\nabla}_k \left\{ g^{kj} \Fvec{j} \right\} + \On{1}{2}  \\
&=& - \frac{1}{2} g^{ij} \Fvec{i} \Fvec{j} -  \nabla_k \left\{ g^{kj}  \Fvec{j}   \right\} + \On{1}{2}  \\
&=& - \half g^{ij} \Fvec{i} \Fvec{j} \\
&&  - \frac{\pi_J}{f} \Delta \frac{f}{\pi_J} + g^{ij} \Fvec{i}\Fvec{j}    + \On{1}{2}  \\
&=& \half g^{ij} \Fvec{i} \Fvec{j} - \frac{\pi_J}{f} \Delta \frac{f}{\pi_J}+ \On{1}{2}.  \\
\end{eqnarray*}
In the above calculation, we used some formulas with respect to the Laplace-Beltrami operator (see Appendix D.).
\\
Thus, we obtain the following theorem, \\
\\
\theorem{1}
 For the ARMA models, if there exists a superharmonic function  $h(\theta) $ such that
$ \Delta h \leq 0$ and $h > 0$, 
then up to the second order, one can improve $\hat{S}_{\pi_J}$ based on the Jeffreys prior 
by adopting the superharmonic prior $\pi_H (\theta) := \pi_{J}(\theta)  h(\theta)$.

\section{Summary}

In the present paper we obtain the asymptotic expansion of the risk difference of the KL-divergence in the ARMA model.
If there exists a superharmonic function  $h(\theta) $ on the corresponding ARMA model manifold, it is better in the Bayesian framework to adopt a superharmonic prior
$\pi_{H}(\theta) := \pi_{J}(\theta) h(\theta)$ as a noninformative prior. It is because that Bayesian spectral densities based on a superharmonic 
prior asymptotically dominates those based on the Jeffreys prior in evaluating the averaged Kullback-Leibler loss.

It is shown that there exists a superharmonic prior for the AR($2$) process and the MA($2$) process~\cite{FT}.
The explicit form of the superharmonic prior is also obtained and the numerical simulation ensures our theorem~\cite{FT2005}.
The existence of superharmonic priors for the higher order ARMA($p$,$q$) processes ($p+q \geq 3$) remains to be discussed.

\vspace{1cm}
\noindent\textbf{Acknowledgment}\\
F.T. was supported by the JSPS Research Fellowships for Young Scientists.

\appendix

\section{Asymptotic Expansion of $\delta $}

In the present section, we evaluate $E_{\theta_0}[\delta^k] $ up to $O_p(n^{-1})$.The key equation is as follows.
\[
 0 = \frac{\partial l (\mle) }{ \partial \theta^m } = 
   \frac{\partial l (\theta_0)}{ \partial \theta^{m} }     
+ \frac{\partial^2 l (\theta_0)}{ \partial \theta^{i} \partial \theta^{m} } \delta^i 
+ \frac{\partial^3 l (\theta_0)}{ \partial \theta^{i} \partial \theta^{j}  \partial \theta^{m} } \frac{1}{2! } \delta^i \delta^j
 + \cdots 
\]
We set $L_m:= \frac{1}{n} \frac{\partial l (\theta_0)}{ \partial \theta^{m} } $ etc.
We omit $\theta_0$ in the remainder of this section.
Since $\delta^i = O_p(n^{-\frac{1}{2} } )   $, higher order terms are recursively obtained and

\[
 \delta^i = - L^{im} L_m - \frac{1}{2} L^{lk}L_{kij} (L^{im}L_m  )  (L^{jn}L_n)   + O_p(n^{- \frac{3}{2}}),
\]
where  $- L^{im} L_m  = O_p(n^{-\frac{1}{2} } ) $ and $ L^{lk}L_{kij} (L^{im}L_m  )  (L^{jn}L_n) =  O_p(n^{-1} ) $.
Now we evaluate $E[\delta^i]$ up to $O_p(n^{-1})$.
Since $ L^{lk}L_{kij} (L^{im}L_m  )  (L^{jn}L_n) =  O_p(n^{-1} ) $, 
\begin{eqnarray*}
&& E[ L^{lk}L_{kij}   (L^{im}  L_m  )  (L^{jn} L_n)] \\
&=&  E[L^{lk}] E[L_{kij} ] E[L^{im} ]E[L^{jn}] E[L_m L_n] + \On{3}{2}
\end{eqnarray*}
Note that some identities $E[L_m] = 0, \quad  E[L_m L_n] = - \frac{1}{n}E[L_{mn}] $.
Let us denote  $m_{ij} := E[L_{ij}] $, $L^{ij} = E[L^{ij}] + \Op{1}{2}$ etc. and $m^{ij}$ be the inverse matrix of $m_{ij}$.
Note that $E[L^{ij}] = m^{ij}+\On{1}{2}$.
Thus,
\begin{eqnarray*}
&&  E[- \frac{1}{2}   L^{lk}L_{kij}   (L^{im}  L_m  )  (L^{jn} L_n)] \\
&=& -\frac{1}{2} E[L^{lk}] E[L_{kij} ] E[L^{im} ]E[L^{jn}] E[L_m L_n] + \On{3}{2}\\
	&=&  -\frac{1}{2n} m^{lk} m_{kij}  m^{im} m^{jn}( -m_{mn}) + \On{3}{2} \\
 	&=&  +\frac{1}{2n} m^{lk} m_{kij} m^{ij} + \On{3}{2}.
\end{eqnarray*}
Collecting the all terms, we can rewrite $E[\delta^{i}]$,
\begin{eqnarray*}
E[\delta^i] &=& E[- L^{im} L_m - \frac{1}{2} L^{lk}L_{kij} (L^{im}L_m  )  (L^{jn}L_n) ]  + O(n^{- \frac{3}{2}}) \\
	&=& -E[L^{im}L_m] + \frac{1}{2n} m^{lk} m_{kij} m^{ij} + \On{3}{2}.
\end{eqnarray*}
Here, the first term is $\On{1}{2}$ and we need to evaluate the second principal term.
\begin{eqnarray*}
E[L^{im}L_m] &=& E\left[ \left\{ m^{im} - m^{il} \delta L_{lk} m^{km} + O_p(n^{-1}) \right\} L_m \right] \\
		&=& m^{im} E[L_m] -  m^{il} E[\delta L_{lk} L_m] m^{km} + \On{3}{2} \\
		&=& -  m^{il} E[\delta L_{lk} L_m] m^{km} + \On{3}{2}. \\
		&=& -  \frac{1}{n} m^{il} m_{lk,m} m^{km} + \On{3}{2},
\end{eqnarray*}
where $ m_{lk,m} := n E[\delta L_{lk} L_m]$.
Thus, we obtain
\[
E[\delta^i] =\frac{1}{n} m^{il} m_{lk,m} m^{km} + \frac{1}{2n} m^{lk} m_{kij} m^{ij} + \On{3}{2}.
\]
\\


\section{Explicit Form of the Expectation of the Log Likelihood for the ARMA model}

In this section, we calculate $m_{ij}$, $m_{ijk}$ and $m_{ij,k}$ for the ARMA model.
Up to $O(n^{-1})$, it can be written in geometrical quantities $g_{ij}$,$T_{ijk}$ and $\mGamma_{i,jk}$, which are defined by spectral density $S(\omega | \theta)$.

\subsection{Trace formula}

Before going into details, we mention the trace formulas. 
Suppose that $\{ x_i \}_{i=1} $ subject to a stationary Gaussian process with zero mean, i.e.,\\  $(x_1, \dots , x_n )\sim \mathcal{N}(0, \Sigma)$.
(The $(s,t)$th component $\Sigma_{st} $ depends only on $s-t$ and such a matrix is called a Toeplitz matrix.)
Then for any symmetric matrices $A$ and $B$, the following equations hold.
\begin{eqnarray*}
E[X_i A_{ij} X_j ] &=& \Tr[ \Sigma A], \\
E[X_i A_{ij} X_j  X_k B_{kl }X_l] &=& \Tr[ A \Sigma ] \Tr[B \Sigma] + 2 \Tr[ A \Sigma B \Sigma ]. 
\end{eqnarray*}

\subsection{$m_{ij}:=E_{\theta_0}[ L_{ij} ] $ for the ARMA model}

In this subsection, we calculate $m_{ij}$.
\begin{eqnarray*}
\frac{\partial^2 l_n}{\partial \theta^i \partial \theta^j  } 
&\!\!\!=&\!\!\! \half X'_n \!\left\{\!\Si \Spp{i}{j} \Si\! \right.  \\
& & \quad \quad - \left. \!\left( \!\Si \Sp{i} \Si \Sp{j} \Si \!+\!\Si \Sp{j} \Si \Sp{i} \Si \right) \!\!  \right\}\! X_n \\
&  &\! + \half  { \Tr \left(\Si \Sp{i} \Si \Sp{j} \right)  - \Tr \left( \Si \Spp{i}{j} \right)   } 
\end{eqnarray*}
We set
\[
\!\!\!S_{ij} :\!=\!\! \frac{1}{2n} \!\!\left\{ \! \Si \Spp{i}{j} \Si \!- \!\left( \!\Si \Sp{i} \Si \Sp{j} \Si \!+\!\Si \Sp{j} \Si \Sp{i} \Si \!\right) \!\!  \right\}, 
\]
\[
 J'_{ij} := \frac{1}{2n}  \Tr \left( \Si \Sp{i} \Si \Sp{j} \right)  \ \mbox{and} \  h_{ij} := \frac{1}{2n} \Tr \left(  \Si \Spp{i}{j} \right).  
\]
Taking expectation of
  $L_{ij} ( :=  \frac{1}{n} \frac{\partial^2 l_n}{\partial \theta^i \partial \theta^j})  = X'_n S_{ij} X_n + J'_{ij} - h_{ij}$,
we obtain
\begin{eqnarray*}
 m_{ij} &:=&  E_{\theta_0}[L_{ij}] \\  
	&=& E[ X'_n S_{ij} X_n + J'_{ij} - h_{ij} ] \\
	&=& \Tr S_{ij}\Sigma + (J'_{ij} - h_{ij}).
\end{eqnarray*}
The first term is rewritten
\begin{eqnarray*}
&& \Tr S_{ij}\Sigma \\
&=& \frac{1}{2n} \Tr \!\left\{  \Si \Spp{i}{j}  - \! \left( \Si \Sp{i} \Si \Sp{j} +\Si \Sp{j} \Si \Sp{i} \right) \!  \right\}\\
&=& h_{ij} - 2 J'_{ij}.
\end{eqnarray*}
Thus,
\begin{eqnarray*}
m_{ij} &=& (h_{ij} - 2 J'_{ij}) + (J'_{ij} - h_{ij}) \\
	&=& - J'_{ij}.
\end{eqnarray*}

\subsection{$m_{ijk}:=E_{\theta_0}[ L_{ijk} ] $ for the ARMA model}

In this subsection, we calculate $m_{ijk}$.
Using the notation in the previous section,  $L_{kij}$ is rewritten by
\[
L_{kij} := \frac{1}{n} \frac{\partial^3 l_n}{\partial \theta^k \partial \theta^i \partial \theta^j  } 
= \frac{\partial }{\partial \theta^k} (X_n' S_{ij} X_n + J'_{ij} - h_{ij}).
\]
Putting the second term and the third term together, we obtain 
\begin{eqnarray*}
&&\frac{\partial }{\partial \theta^k}J'_{ij} - \frac{\partial }{\partial \theta^k} h_{ij} \\
 &=& \frac{1}{2n}  \Tr \left( \Si \Spp{k}{i} \Si \Sp{j} + \Si \Sp{i} \Si \Spp{k}{j} \right)  \\
 && - \frac{1}{2n} \Tr \left( \Si \Sp{k} \Si \Sp{i} \Si \Sp{j} \Si + \Si \Sp{i} \Si \Sp{k} \Si \Sp{j} \Si  \right) \\
&&  - \frac{1}{2n} \Tr \left( \Si \Sppp{k}{i}{j} \right) +\frac{1}{2n} \Tr \left( \Si \Sp{k} \Si \Spp{i}{j} \right) \\
 &=& (\Gamma'_{jki} + \Gamma'_{ikj} + \Gamma'_{kij} )  - (T'_{kij} + T'_{ikj}) - N'_{kij},
\end{eqnarray*}
where
\[
 \Gamma'_{kij} := \Tr\left( \Si \Sp{k} \Si \Spp{i}{j} \right), 
\]
\[
 T'_{kij} := \Tr \left(\Si \Sp{k} \Si \Sp{i} \Si \Sp{j} \right),
\]
and
\[
N'_{kij} := \Tr\left(\Si \Sppp{k}{i}{j} \right).
\]
Since the first term is in the tedious form, we introduce the following notation:
\begin{eqnarray*}
&&A_{ijk} + \mbox{(permutation terms)} := A_{ijk} + A_{ikj} + A_{jki}+A_{jik} + A_{kij}+A_{kji}, \\
&&B_{ijk} +  \mbox{(cyclic terms)} := B_{ijk} + B_{jki} + B_{kij}.
\end{eqnarray*}
Then, 
\begin{eqnarray*}
S_{kij}
 &:=& \frac{\partial}{ \partial \theta^k} S_{ij} \\
&=&\frac{1}{2n} \left\{  \Si \Sp{k} \Si \Sp{i} \Si \Sp{j} \Si + \mbox{(permutation terms)} \right\} \\
 & & - \frac{1}{2n} \left\{ \left(  \Si \Spp{k}{i}  \Si \Sp{j} \Si +\Si \Sp{i}  \Si \Spp{k}{j} \Si \right) \right. \\
 &&	\quad \quad \quad + \mbox{(cyclic terms)}  \Bigg\}  + \frac{1}{2n} \left(\Si \Sppp{k}{i}{j} \Si \right).
\end{eqnarray*}
Taking average, we obtain
\begin{eqnarray*}
&&  E[X'_n S_{kij} X_n] \\
&=& \Tr S_{kij} \Sigma \\
&=&\frac{1}{2n} \Tr \left\{  \left( \Si \Sp{i}\Si \Sp{j} \Si \Sp{k} \right) + \mbox{(permutation terms)}  \right\}  \\
& & -\frac{1}{2n}  \Tr \left\{  \left( \Si \Spp{k}{i}  \Si \Sp{j} +\Si \Sp{i}  \Si \Spp{k}{j} \right) \right. \\
& & \quad \quad  \quad \ \ + \mbox{(cyclic terms)} \Bigg\}+ \frac{1}{2n} \Tr \left( \Si \Sppp{k}{i}{j} \right) \\
&=& T'_{ijk}+T'_{jki}+T'_{kij} + T'_{jik}+T'_{ikj}+T'_{kji}
	-2 ( \Gamma'_{ijk} + \Gamma'_{jik} + \Gamma'_{kij} ) + N'_{kij} \\
&=& 3(T'_{ijk}+T'_{jik}) -2 ( \Gamma'_{ijk} + \Gamma'_{jik} + \Gamma'_{kij} ) + N'_{ijk}.
\end{eqnarray*}
In the last line, the cyclic property $T'_{ijk} = T'_{jki}=T'_{kij}$, etc. was used, which is due to the property of the trace operation $\Tr ABC = \Tr BCA = \Tr CAB $. 
\\
\\
Thus, $m_{kij}$ is written in the form of
\begin{eqnarray*}
m_{kij} &:=& E[L_{kij}] \\
 &=&  E[X'_n S_{kij} X_n] +  (\Gamma'_{kij} + \Gamma'_{jki} + \Gamma'_{ikj} )
 - (T'_{kij} + T'_{ikj}) - N'_{ijk} \\
 &=&  3(T'_{ijk}+T'_{jik}) -2 ( \Gamma'_{ijk} + \Gamma'_{jik} + \Gamma'_{kij} ) + N'_{ijk} \\
 &&  +  (\Gamma'_{kij} + \Gamma'_{jki} + \Gamma'_{ikj} ) - (T'_{kij} + T'_{ikj}) - N'_{ijk} \\
&=&  2(T'_{ijk}+T'_{jik}) - ( \Gamma'_{ijk} + \Gamma'_{jik} + \Gamma'_{kij} ).
\end{eqnarray*}

\subsection{$m_{ij,k}:=E_{\theta_0}[ L_{ij} L_{k}] $ for the ARMA model}

In this section we calculate $m_{ij,k}:=E_{\theta_0}[ L_{ij} L_{k}]$.
\begin{eqnarray*}
L_{ij}  &=&  X'_n S_{ij} X_n + J'_{ij} - h_{ij}\\
L_k  &=& \frac{1}{2n} \left\{ X'_n \left( \Sigma^{-1} \Sp{k} \Sigma^{-1} \right) X_n - \Tr \left(\Sigma^{-1} \Sp{k}\right) \right\} \\
       &=& X'_n B_k X_n - E[X'_n B_k X_n],
\end{eqnarray*}
where
\[
 B_k :=  \frac{1}{2n} \left\{ \Sigma^{-1} \left(\Sp{k} \right)\Sigma^{-1} \right\}.
\]
Thus, 
\begin{eqnarray*}
&&m_{ij,k} \\
&=& E[ (X'_n S_{ij} X_n + J'_{ij} - h_{ij})  ( X'_n B_k X_n - \Tr [ B_k \Sigma] )  ]\\
	&=&\! \Tr[ S_{ij}\Sigma] \Tr[B_k \Sigma] + 2 \Tr[ S_{ij}\Sigma B_k \Sigma]  - \Tr[S_{ij}\Sigma] \Tr[B_k \Sigma] \\ 
	&=&\!  2 \Tr[ S_{ij}\Sigma B_k \Sigma] \\
	&=&\!  2\Tr\left[   \frac{1}{2n} \left\{  \Si \Spp{i}{j} - \left( \Si \Sp{i} \Si \Sp{j} \Si +\Si \Sp{j} \Si \Sp{i} \right)   \right\} \right.
             \\ && \quad \times  \left. \frac{1}{2n} \left\{ \Si \Sp{k} \right\} \right]    \\
	&=&\! \frac{1}{n} \left\{   \frac{1}{2n}\Tr\left( \Si \Spp{i}{j} \Si \Sp{k} \right) -\frac{1}{2n}\Tr\left( \Si \Sp{i} \Si \Sp{j} \Si \Sp{k} \right)  \right. \\
	&&\!  {} \left. - \frac{1}{2n}\Tr\left( \Si \Sp{j} \Si \Sp{i} \Si \Sp{k} \right)\right\} \\
	&=&\! \frac{1}{n} ( \Gamma'_{kij} - T'_{ijk} - T'_{jik} )  \ ( = O(n^{-1}) ) 
\end{eqnarray*} 
Now we summarize the whole results,
\[
\begin{array}{ccl}
 m_{ij}   &=& - J'_{ij},  \\ 
 m_{ijk}  &=& 2(T'_{ijk}+T'_{jik}) - ( \Gamma' _{ijk} +\Gamma' _{jik} + \Gamma'_{kij} ), \\ 
 m_{ij,k} &=& \frac{1}{n} (\Gamma'_{kij} - T'_{ijk} - T'_{jik}). \\
\end{array}
\]


\section{Moment Formula}

Let us calculate the $p$-th moment of the multivariate Gaussian distribution by using the characteristic function (chf).
We set a mean parameter equal to zero and denote a covariance matrix as $\Sigma $.
Then, chf is given by
\[
\varphi(t) := E[e^{it'X}] = \exp\left(-\frac{1}{2} t' \Sigma t\right)
\]
where $t = (t^1 \cdots t^k )' $.For even $p$,
\[
 I^{i_1 i_2 \cdots i_p} = (i)^p  \frac{\partial }{\partial t^{i_1} }  \frac{\partial }{\partial t^{i_2} }
  \cdots \frac{\partial }{\partial t^{i_p} }  
 \left.
	\exp\left(-\frac{1}{2} t' \Sigma t\right)
 \right|_{t \rightarrow 0}
\]
Note that for odd $p$, the moments vanish.

\section{Laplace-Beltrami Operator}

We briefly summarize the Laplace-Beltrami operator in the Riemannian manifold. (See, e.g., Kobayashi and Nomizu~\cite{KN}).
The covariant derivative in the $j$-th direction of a vector $V^{l}$ is defined by
\[
\nabla_j V^l := \partial_j V^l + \Gamma^{l}_{kj}V^k, 
\]
where $\partial_j $ denotes $\frac{\partial}{\partial \theta^{j}}$.
When we set $V^j = \nabla^j \phi = \partial^j \phi = g^{jl} \partial_{l} \phi $ for a scalar function $\phi $,
the Laplace-Beltrami operator is defined by 
\begin{eqnarray*}
\Delta \phi  := \nabla_{j} \nabla^j \phi &=& \partial_j (\nabla^j \phi ) + \Gamma_{jk}^k (\nabla^j \phi ) \\
		&=& \partial_j (\nabla^j \phi ) + \frac{1}{\sqrt{g}} (\partial_j \sqrt{g}) (\nabla^j \phi ) \\
		&=& \frac{1}{\sqrt{g}} \partial_i \left(  \sqrt{g}  \nabla^j \phi  \right) \\
		&=& \frac{1}{\sqrt{g}} \partial_i \left(  \sqrt{g}  g^{jk} \partial_{k} \phi     \right). 
\end{eqnarray*}
Denoting $g := \det{g_{ij}}$, and since $\frac{\partial_j g}{g} = g^{kl} \partial_j g_{kl}$,
\[
 \Gamma^i _{ji} = \half g^{ik} (\partial_i g_{kj} + \partial_j g_{ki} - \partial_k g_{ji})
	= \half g^{ik} \partial_jg_{ki}
	= \half \frac{\partial_j g}{ g}.
\]
Thus we can rewrite $\Gamma^{i}_{ji} = \frac{1}{\sqrt{g}} \partial_j (\sqrt{g}) = \partial_j \log(\sqrt{g}) = \partial_j \log \pi_J $.

\end{document}